\documentclass[twoside,bezier,reqno]{amsart}

\usepackage{epsfig}
\usepackage{amsmath,amsthm,amsfonts,amssymb,amscd,euscript}
\usepackage{enumerate}
\input{epsf}
\newcommand{\filebegin}{\begin{document}}
\newcommand{\fileend}{\end{document}}
\makeatother
\def\thefootnote{}
\newcommand{\lo}{\longrightarrow}
\newcommand{\NMM}{\hspace*{2mm}}

\renewcommand{\baselinestretch}{1.1}

\input{epsf}
\renewcommand{\baselinestretch}{1.1}
\def\n{\noindent}%

\numberwithin{equation}{section}
\def\mapdown#1{\Big\downarrow\rlap
{$\vcenter{\hbox{$\scriptstyle#1$}}$}}
\newtheorem{theorem}{Theorem}[section]
\newtheorem{lemma}[theorem]{Lemma}
\newtheorem{proposition}[theorem]{Proposition}
\newtheorem{corollary}[theorem]{Corollary}
\newtheorem{conjecture}{Conjecture}[section]
\newtheorem{problem}{Problem}[section]

\theoremstyle{definition}
\newtheorem{definition}[theorem]{Definition}
\newtheorem{example}[theorem]{\sc Example}
\newtheorem{xca}[theorem]{Exercise}
\theoremstyle{remark}
\newtheorem{remark}[theorem]{Remark}
\usepackage{xcolor}
\definecolor{vividviolet}{rgb}{0.62, 0.0, 1.0}
\def\Shiu{\color{vividviolet}}
\def\rsq{\hspace*{\fill}$\blacksquare$}

\begin{document}

\setcounter{page}{1} \noindent
Iranian Journal of Mathematical Sciences and Informatics \\
Vol. x, No. x (202x), pp xx-xx

\vspace*{2cm}
\begin{center}
{\bf\large On local antimagic chromatic number of graphs with cut-vertices}
 \\[0.5cm]
{Gee-Choon Lau$^a*$, Wai-Chee Shiu$^{b}$, Ho-Kuen Ng$^{c}$ \footnote{$^*$Corresponding Author} \\[2mm]
$^a$Faculty of Computer \& Mathematical Sciences, Universiti Teknologi MARA (Segamat Campus), 85000 Johor, Malaysia\\
$^b$Department of Mathematics, The Chinese University of Hong Kong, Shatin, Hong Kong\\
$^c$Department of Mathematics, San Jos\'{e} State University, San Jos\'{e} CA 95192 USA\\[2mm]
{\tt E-mail: geeclau@yahoo.com}\\
{\tt E-mail: wcshiu@associate.hkbu.edu.hk}\\
{\tt E-mail: ho-kuen.ng@sjsu.edu}
} \\[2mm]
\end{center}%
\vspace*{0.5cm}
\begin{quotation}
\noindent
{\footnotesize
{\sc Abstract.}
An edge labeling of a connected graph $G = (V, E)$ is said to be local antimagic if it is a bijection $f:E \to\{1,\ldots ,|E|\}$ such that for any pair of adjacent vertices $x$ and $y$, $f^+(x)\not= f^+(y)$, where the induced vertex label $f^+(x)= \sum f(e)$, with $e$ ranging over all the edges incident to $x$.  The local antimagic chromatic number of $G$, denoted by $\chi_{la}(G)$, is the minimum number of distinct induced vertex labels over all local antimagic labelings of $G$. In this paper, the sharp lower bound of the local antimagic chromatic number of a graph with cut-vertices given by pendants is obtained. The exact value of the local antimagic chromatic number of many families of graphs with cut-vertices (possibly given by pendant edges) are also determined. Consequently, we partially answered Problem 3.1 in [Local antimagic vertex coloring of a graph, {\it Graphs and Combin.}, {\bf33} (2017),  275--285.].}
\end{quotation}

\n \textbf{2000 Mathematics subject classification: } 05C78, 05C69.

\markboth
{Gee-Choon Lau, Wai-Chee Shiu, Ho-Kuen Ng}
 {On local antimagic chromatic number of graphs with cut-vertices}



\section{Introduction}

 A connected graph $G = (V, E)$ is said to be {\it local antimagic} if it admits a {\it local antimagic edge labeling}, i.e., a bijection $f : E \to \{1,\dots ,|E|\}$ such that the induced vertex labeling $f^+ : V \to \mathbb{Z}$ given by $f^+(u) = \sum f(e)$ (with $e$ ranging over all the edges incident to $u$) has the property that any two adjacent vertices have distinct induced vertex labels. Thus, $f^+$ is a coloring of $G$. Clearly, the order of $G$ must be at least 3.  The vertex label $f^+(u)$ is called the {\it induced color} of $u$ under $f$ (the {\it color} of $u$, for short, if no ambiguous occurs). The number of distinct induced colors under $f$ is denoted by $c(f)$, and is called the {\it color number} of $f$. The {\it local antimagic chromatic number} of $G$, denoted by $\chi_{la}(G)$, is $\min\{c(f) : f\mbox{ is a local antimagic labeling of } G\}$. Clearly, $2\le \chi_{la}(G)\le |V(G)|$. The sharp lower bound of the local antimagic chromatic number of a graph with cut-vertices given by pendants is obtained. In~\cite[Problem 3.3]{Arumugam}, the authors asked:

\begin{center} Does there exist a graph $G$ of order $n$ with $\chi_{la}(G) = n - k$ for every
$k = 0, 1, 2, \ldots, n - 2$? \end{center}

   In~\cite[Theorems 3.4 and 3.5]{LNS}, we proved the following that answered the above problem affirmatively.

\begin{theorem}{\rm\cite{LNS}} For each possible $n,k,$ there exists a graph $G$ of order $n$ such that $\chi_{la}(G) = n - k$ if and only if $n\ge k+3\ge 3$. Moreover, there is a graph $G$ of order $n$ with $\chi_{la}(G)=2$ if and only if $n\ne 2, 3, 4, 5, 7$. \end{theorem}

 We shall in Section~2 completely determine the local antimagic chromatic number of the one-point union of cycles. Let $G$ be a graph of order $n\ge 3$. We also determined the exact value of the local antimagic chromatic number of many families of graphs with pendants that has $\chi_{la}(G)<n$. In Section~3, we obtained several families of graphs $G$ with $\chi_{la}(G)=n$. This partially answered~\cite[Problem 3.1]{Arumugam}. For convenience, we shall use $a^{[n]}$ to denote a sequence of length $n$ in which all items are $a$, where $n\ge 2$. For integers $1\le a < b$, we let $[a,b]$ denote the set of integers from $a$ to $b$.

\section{$\chi_{la}(G) < |V(G)|$}

 In~\cite{Arumugam}, the authors proved that for every tree $T$ with $k$ pendant edges (i.e., with $k$ pendants), $\chi_{la}(T)\ge k+1$. We generalize this result to arbitrary graphs of order at least 3.

\begin{lemma}\label{lem-pendant} Let $G$ be a graph of size $q$ containing $k$ pendants. Let $f$ be a local antimagic labeling of $G$ such that $f(e)=q$. If $e$ is not a pendant edge, then $c(f)\ge k+2$.\end{lemma}
\begin{proof}
Let $e=uv$ and $x_1, \dots, x_k$ be pendants. Thus, $f^+(u)>q$ and $f^+(v)>q$ and they are distinct. On the other hand, $f^+(x_i)<q$ and are distinct for all $i$. Hence $c(f)\ge k+2$.
\end{proof}

\begin{theorem}\label{thm-pendant}  Let $G$ be a graph having $k$ pendants. If $G$ is not $K_2$, then $\chi_{la}(G)\ge k+1$ and the bound is sharp.\end{theorem}

\begin{proof} Suppose $G$ has size $q$. Let $f$ be any local antimagic labeling of $G$.
Consider the edge $uv$ with $f(uv) = q$. We may assume $u$ is not a pendant. Clearly, $f^+(u) > q \ge f^+(z)$ for every pendant $z$. Since all pendants have distinct induced colors, we have $\chi_{la}(G)\ge k+1$.

 For $k\ge 2$, since $\chi_{la}(S_k)= k+1$, where $S_k$ is a star with maximum degree $k$, the lower bound is sharp. The left labeling below is another example also showing that the lower bound is sharp. The right labeling shows that the lower bound is sharp for $k=1$.\\\\
\centerline{\epsfig{file=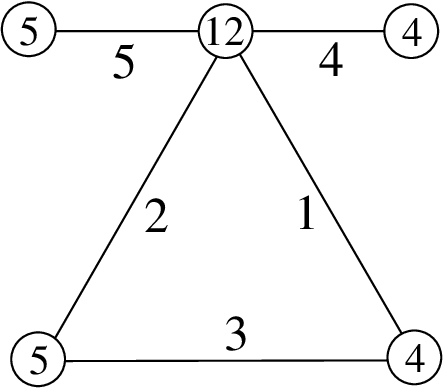, width=2.5cm}\qquad \epsfig{file=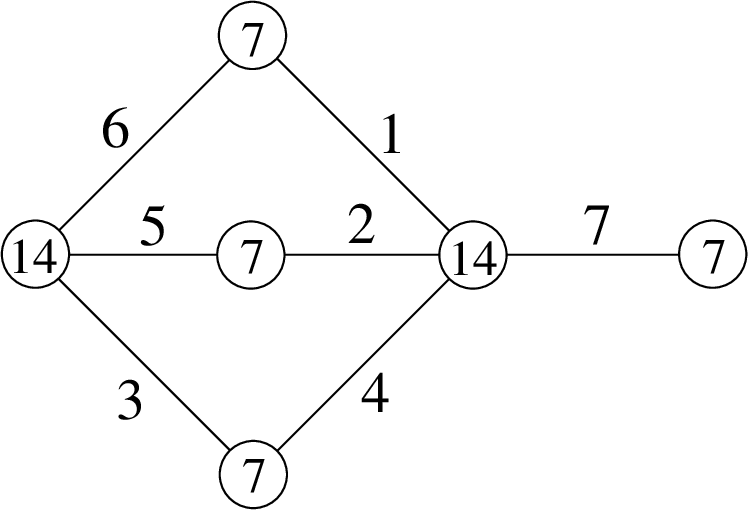, width=3.3cm}}
 \end{proof}

 The contrapositive of the following lemma \cite[Lemma~1]{LNS} or \cite[Lemma~2.1]{LSN} gives a sufficient condition for a bipartite graph $G$ to have $\chi_{la}(G)\ge 3$.

\begin{lemma}{\rm\cite{LNS,LSN}}\label{lem-2part} Let $G$ be a graph of size $q$. Suppose there is a local antimagic labeling of $G$ inducing a $2$-coloring of $G$ with colors $x$ and $y$, where $x<y$. Let $X$ and $Y$ be the sets of vertices colored $x$ and $y$, respectively, then $G$ is a bipartite graph with bipartition $(X,Y)$ and $|X|>|Y|$. Moreover,
$x|X|=y|Y|= \frac{q(q+1)}{2}$.
\end{lemma}

 For $r\ge 2$ and $a_1\ge a_2\ge \cdots\ge a_r\ge 3$, denote by $C(a_1,a_2,\ldots,a_r)$ the one-point union of $r$ distinct cycles of order $a_1, a_2,\ldots,a_r$ respectively. Note that $C(a_1,a_2,\ldots,a_r)$ has $m=a_1 + \cdots + a_r\ge 6$ edges and $m-r+1$ vertices. We shall denote the vertex of maximum degree by $u$, called the {\it central vertex}, and the $2r$ edges incident to $u$ are called the {\it central edges}. Denote the consecutive edges of subgraph $C_{a_i}$ by $e_{s_i+1}, e_{s_i+2}, \ldots, e_{s_i+a_i}$ such that $s_1=0$, $s_i=a_1+a_2+\cdots+a_{i-1}$ for $i\ge 2$. Moreover, for $i\ge 1$, $e_{s_i+1}$ and $e_{s_i+a_i}$ are the central edges of $C_{a_i}$.

\begin{theorem}\label{thm-vertexgluecycles}  Suppose $G=C(a_1,a_2,\ldots,a_r)$, then $\chi_{la}(G)=  2$ if and only if $G=C((4r-2)^{[r-1]},2r-2)$, $r\ge 3$ or $G=C((2r)^{[(r-1)/2]}$, $(2r-2)^{[(r+1)/2]})$, $r$ is odd. Otherwise, $\chi_{la}(G)= 3$. \end{theorem}

\begin{proof} Let $G=C(a_1,a_2,\ldots,a_r)$. Define an edge labeling $f: E(G)\to [1,m]$ by
\begin{enumerate}[1.]
  \item $f(e_i) = i/2$ for even $i$,
  \item $f(e_i) = m - (i-1)/2$ for odd $i$.
\end{enumerate}
 It is easy to verify that $f^+(u)> m+1$, and each vertex of degree 2 has color $m+1$ and $m$ alternately beginning from vertices adjacent to $u$. Therefore, $f$ is a local antimagic labeling that induces a 3-coloring. Thus, $\chi_{la}(G)\le 3$. If $G$ contains an odd cycle, we have $\chi_{la}(G)\ge\chi(G)= 3$ so that $\chi_{la}(G)=3$.

Suppose $\chi_{la}(G)=2$. This implies that $\chi(G)=2$ and hence $a_i\ge 4$ is even for each $i$. Let $g$ be any local antimagic coloring of $G$ that induces a 2-coloring of $G$ with colors $x$ and $y$. Without loss of generality, we may assume that $g^+(u)=y$. Let $X$ and $Y$ be the sets of vertices with colors $x$ and $y$, respectively. It is easy to get that $|Y|=m/2-r+1$ and $|X|=m/2$. By Lemma~\ref{lem-2part}, we have $x|X|=y|Y|=m(m+1)/2$. Hence, $x=m+1\ge 4r+1$ is odd,  $y=m(m+1)/(m-2r+2)$ and $y\ge 1+2+\cdots+2r= 2r^2+r$. Suppose $\ell$ is labeled at an edge $vw$ which is not a central edge. Without loss of generality, we may assume that $g^+(v)=x$ and $g^+(w)=y$. Then the label assigned to another edge incident with $w$ must be $y-\ell$. Then $1\le y-\ell\le m$, i.e., $\ell\ge y-m=y-x+1$. In other word, labels in $[1, y-x]$ are labeled at central edges. So $y-x\le 2r$.

\medskip

Solving for $m$, we get $m = (y-1\pm \sqrt{y^2+6y-8yr+1})/2$. Hence, $y^2+6y-8yr+1 = t^2 \ge 0$, where $t$ is a nonnegative integer. This gives $(y+3-4r)^2 + 1 - (3-4r)^2 = t^2$ or $(y+3-4r-t)(y+3-4r+t) = 8(2r-1)(r-1)$. By letting $a=y+3-4r-t$ and $b=y+3-4r+t$ we have $2y+6-8r=a+b$ with $ab=8(r-1)(2r-1)=8(2r^2-3r+1)$. Clearly $b\ge a>0$. Since $a,b$ must be of same parity, we have both $a,b$ are even.

Recall that $y-2r^2-r\ge 0$. Now
\begin{align*}y-2r^2-r& =4r-3+\frac{a+b}{2}-2r^2-r\\
& = \frac{a+b}{2} -2r^2+3r -3=\frac{a+b}{2} -\frac{ab}{8}-2\\
& = \frac{4a+4b-ab-16}{8}= -\frac{(a-4)(b-4)}{8}.
\end{align*}
This implies that $a\le 4$.

 Before considering the cases when $a\le 4$, we need the following claim which is easy to obtain.\\
\noindent{\bf Claim: }{\it Let $\phi$ be a labeling of a $2s$-cycle $v_1v_2\cdots v_{2s}v_1$ with $\phi(v_{2i-1}v_{2i})=\alpha_i$ and $\phi(v_{2i}v_{2i+1})=\beta_i$ for $1\le i\le s$, where $v_{2s+1}=v_1$. Suppose $\phi^+(v_{2j})=x$ for $1\le j\le s$ and $\phi^+(v_{2k+1})=y$ for $1\le k\le s-1$, where $y>x$. Then $\alpha_1+\beta_1=x$, $\{\alpha_1,\alpha_2,\dots,\alpha_s\}$ is an increasing sequence with common difference $y-x$ and $\{\beta_1,\beta_2,\dots,\beta_{s}\}$ is an decreasing sequence with common difference $y-x$.}

\noindent{\bf Case (1).} Suppose $a=2$. In this case, $b=4(r-1)(2r-1)$ and $2y+6-8r=8r^2-12r+6$. Hence, $y=4r^2-2r$. This gives (i) $m=4r^2-4r$ and $x=4r^2-4r+1$ or (ii) $m=2r-1$ and $x=2r < 4r+1$, a contradiction. In (i), $y-x=2r-1$. Since $[1,2r-1]$ must be assigned to central edges, the central edges must be labeled by 1 to $2r-1$ and $2r^2-r$, respectively. There are $r-1$ cycles, say $C_{a_1}, \dots, C_{a_{r-1}}$, whose central edges are labeled by numbers in $[1, 2r-1]$.

It is easy to verify  that no such graph exists for $r=2$. So we assume that $r\ge 3$.

Suppose $C_{2s}$ is one of these $r-1$ cycles. Note that $s\ge 2$.  Keep the notation defined in the claim. By symmetry, we may assume that $\alpha_1<\beta_s$. So $\alpha_1\in[1, 2r-2]$. Now we have $\beta_{s}=(x-\alpha_1)-(s-1)(y-x)\le 2r-1$ and $\beta_{s-1}=(x-\alpha_1)-(s-2)(y-x)\ge 2r$. Thus, $(2r-1)^2-\alpha_1\le s(y-x)\le 2r-2+(2r-1)^2-\alpha_1$. Hence, $(2r-1)^2-(2r-2)\le s(2r-1)\le (2r-1)^2+(2r-1)-2$. This implies that $2r-2<s<2r$. Thus $s=2r-1$. Moveover, $\beta_s=\beta_{2r-1}=(2r-1)-\alpha_1\le 2r-2$. So, $a_j=4r-2$ for $1\le j\le r-1$.

We are now left with one unlabeled cycle, also denoted by $C_{2s}$, with central edge labels must be $2r-1$ and $2r^2-r$. Again, we may assume $2r-1=\alpha_1<\beta_s=2r^2-r$. By the claim, $2r^2-r=(x-\alpha_1)-(s-1)(y-x)=(2r-1)^2-(2r-1)-(s-1)(2r-1)$. Thus $s=r-1$ and hence $a_r=2r-2$. Therefore, $G=C((4r-2)^{[r-1]},2r-2)$.

On the other hand, for $i$-th $(4r-2)$-cycle, we choose $\alpha_1=i$, $1\le i\le r-1$; for the $(2r-2)$-cycle, we choose $\alpha_1=2r-1$. Apply the labeling as shown in the claim. One can verify that the edge labels are all distinct in $[1,4r^2-4r]$. Consequently, $C((4r-2)^{[r-1]},2r-2)$ admits a local antimagic labeling that induces a 2-coloring. Thus, $\chi_{la}(C((4r-2)^{[r-1]},2r-2))=2$.

\noindent {\bf Case (2).} Suppose $a=4$. In this case, $b=2(r-1)(2r-1)$ and $2y+6-8r=4r^2-6r+6$. Hence, $y=2r^2+r$. This gives (i) $m=2r^2-r-1$ and $x=2r^2-r$ with $r$ is odd or (ii) $m=2r$ and $x=2r+1 < 4r+1$, a contradiction. In (i), $y-x=2r$. Thus all the central edges must be assigned with integers in $[1,2r]$. Suppose $C_{2s}$ is one of the cycles whose central edges are labeled by $\alpha_1$ and $\beta_s$. Also, by symmetry we may assume $\alpha_1<\beta_s$. So $\alpha_1\in[1,2r-1]$. By a similar computation as in Case~(1), we have $r-\frac{3}{2}+\frac{1}{r}\le s\le r+\frac{1}{2}-\frac{1}{r}$. So $s=r-1$ or $r$. Suppose there are $k$ cycles of $2r$ edges in $G$, then there are $r-k$ cycles of $2r-2$ edges. Now, the size of $G$ is $m=2rk+(2r-2)(r-k)=2r^2-2r+2k$. Thus we have $k=\frac{r-1}{2}$. Hence $G=C((2r)^{[(r-1)/2]},(2r-2)^{[(r+1)/2]})$.

Moreover, when $s=r$ and since $\alpha_1<\beta_s$, we have $\alpha_1\le r/2$. Since $r$ is odd, $\alpha_1\le \frac{r-1}{2}$. Thus labels in $[1,\frac{r-1}{2}]$ are labeled at each of $2r$-cycle, respectively.

On the other hand, for $i$-th $2r$-cycle, we choose $\alpha_1=i$, $1\le i\le \frac{r-1}{2}$; for the $j$-th $(2r-2)$-cycle we choose $\alpha_1=j+\frac{r-1}{2}$, $1\le j\le \frac{r+1}{2}$. Apply the labeling as shown in the claim. One may verify that the edge labels are all distinct in $[1,2r^2-r-1]$.  Consequently, $C((2r)^{[(r-1)/2]},(2r-2)^{[(r+1)/2]})$ admits a local antimagic labeling that induces a 2-coloring. Thus, $$\chi_{la}(C((2r)^{[(r-1)/2]},(2r-2)^{[(r+1)/2]}))=2.$$

Consequently, $\chi_{la}(G)=  2$ if and only if $G=C((4r-2)^{[r-1]},2r-2), r\ge 3$ or $G=C((2r)^{[(r-1)/2]}, (2r-2)^{[(r+1)/2]}), r$ is odd. Otherwise, $\chi_{la}(G)= 3$.
\end{proof}

\begin{example} For $C(10,10,4)$, beginning and ending with central edges, the two 10-cycles has consecutive labels 1, 24, 6, 19, 11, 14, 16, 9, 21, 4 and 2, 23, 7, 18, 12, 13, 17, 8, 22, 3 respectively while the 4-cycle has consecutive labels $5, 20, 10, 15$ with $y=30$ and $x=25$. For $C(14,14,14,6)$, the three 14-cycles has consecutive edge labels 1, 48, 8, 41, 15, 34, 22, 27, 29, 20, 36, 13, 43, 6; 2, 47, 9, 40, 16, 33, 23, 26, 30, 19, 37, 12, 44, 5 and 3, 46, 10, 39, 17, 32, 24, 25, 31, 18, 38, 11, 45, 4 respectively, while the 6-cycle has consecutive edge labels $28,21,35,14,42,7$ with $y=56$, and $x=49$. Similarly, for $C(6,4,4)$, the 6-cycle has consecutive edge labels $1,14,7,8,13,2$ while the two 4-cycles has consecutive edge labels $3,12,9,6$ and $4,11,10,5$ respectively with $y=21$ and $x=15$. For $C(10,10,8,8,8)$, the two 10-cycles has consecutive labels 1, 44, 11, 34, 21, 24, 31, 14, 41, 4 and 2, 43, 12, 33, 22, 23, 32, 13, 42, 3 respectively, while the three 8-cycles has consecutive labels $5,40,15,30,25,20,35,10$; 6, 39, 16, 29, 26, 19, 36, 9 and 7, 38, 17, 28, 27, 18, 37, 8 respectively, with $y=55$ and $x=45$. \rsq
\end{example}

 For $k,r\ge 1$ and $a_1\ge a_2\ge \cdots\ge a_r\ge 3$, let $H(a_1,a_2,\ldots,a_r;k)$ be the {\it hibiscus graph} obtained by identifying the central of $C(a_1,a_2,\ldots,a_r)$ with an end-vertex of $k$ copies of $P_2$. Clearly, $H(a_1,a_2,\ldots,a_r;k)$ has $m+k = a_1 + \cdots + a_r + k\ge 4$ edges and $m+k-r+1$ vertices. For non-pendant vertices and edges, we shall adopt the notation of $C(a_1,a_2,\ldots,a_r)$ accordingly.

\begin{theorem} For $k\ge 1$,
\[\chi_{la}(H(a_1,a_2,\ldots,a_r;k))=\begin{cases} 3 & \mbox{if }k=1,\\
k+1 & \mbox{if }k\ge 2.
\end{cases}\]
\end{theorem}

\begin{proof} Let $v_j$ $(1\le j\le k)$ be the pendant vertices of $G=H(a_1,a_2,\ldots,a_r;k)$.
Define an edge labeling $f: E(G)\to [1,m+k]$ by

\begin{enumerate}[1.]
\item $f(e_i)=(i+1)/2$ for odd $i$,
\item $f(e_i) = m-i/2+1$ for even $i$,
\item $f(uv_j) = m + j$ for $1\le j\le k$.
\end{enumerate}

 It is easy to verify that $f^+(u)>m+k+3$, $f^+(v_j)=m+j$ for $1\le j\le k$, and each degree 2 vertex has color $m+1$ and $m+2$ alternately beginning from vertices adjacent to $u$.

When $k\ge 2$, we have that $f$ is a local antimagic labeling that induces a $(k+1)$-coloring so that $\chi_{la}(G)\le k+1$. By Theorem~\ref{thm-pendant}, we know $\chi_{la}(G)\ge k+1$. Therefore, $\chi_{la}(G)=k+1$.

Suppose $k=1$. Clearly, $f$ is a local antimagic labeling that induces a 3-coloring. So $\chi_{la}(G)\le 3$. If $G$ contains an odd cycle, then $\chi_{la}(G)\ge \chi(G) = 3$.  Hence, $\chi_{la}(G)=3$. Suppose $\chi_{la}(G)=2$. Then $G$ is bipartite and hence $a_i$ is even for each $1\le i\le r$. Let $g$ be a local antimagic coloring of $G$ that induces a 2-coloring with colors $x$ and $y$ such that $g^+(u)=y$. By Lemma~\ref{lem-pendant} $g(uv_1)=m+1$. Since $g^+(u)=y$, $g^+(v_1)=x$. Hence $x=m+1$.

Let $X$ and $Y$ be the sets of vertices with colors $x$ and $y$, respectively. It is easy to get that $|Y|=(m-2r+2)/2$ and $|X|=(m+2)/2$.  By Lemma~\ref{lem-2part}, we have $y(m-2r+2)/2=(m+1)(m+2)/2$. Hence, $y=(m+1)(m+2)/(m-2r+2)$.

Solving for $m$, we get $m = (y-3\pm \sqrt{y^2+2y-8yr+1})/2$. Hence, $y^2+2y-8yr+1 = t^2 \ge 0$. This gives $(y+1-4r)^2 - (1-4r)^2 + 1 = t^2$ or $(y+1-4r-t)(y+1-4r+t)=8r(2r-1)$, where $t\ge 0$. By letting $a=y+1-4r-t$ and $b=y+1-4r+t$, we have $2y+2-8r=a+b$ with $ab=8r(2r-1)$. Since $a,b$ must be of same parity, we have both $a,b$ are even.


Now, $y\ge m+1+\sum\limits_{i=1}^{2r} i\ge 4r+1+r(2r+1)=2r^2+5r+1$. By a similar computation in the proof of Theorem~\ref{thm-vertexgluecycles}, $0<2r\le y-2r^2-3r-1=-\frac{(b-4)(a-4)}{8}$. This implies that $a=2$ and $b\ge 6$. In this case, $b=4r(2r-1)$ and hence $y=4r^2+2r$. Thus $t=4r^2-2r-1$ and hence $m=4r^2-2$ or $m=2r-1$. Since $m\ge 4r$, $m=4r^2-2$. Since $b\ge 6$, $r\ge2$.
Now $y=4r^2+2r\ge (4r^2-1)+\sum\limits_{i=1}^{2r} i =6r^2+r-1$ yields a contradiction. Thus $\chi_{la}(H(a_1,a_2,\ldots,a_r;1))=3$. \end{proof}






 Let $T(m,n)$ be the vertex-gluing of the end vertex of a path $P_m$ and a vertex of a cycle $C_n$. In some article, $T(m,n)$ is called a {\it tadpole graph}.

\begin{theorem} For $n \ge 3$, $m\ge 2$, $\chi_{la}(T(m,n)) = 3$.  \end{theorem}

\begin{proof} Note that $T(m,n)$ has order and size $m+n-1$. Let the edge set be $\{e_i=v_iv_{i+1} \;|\; i\in[1,m+n-2]\}\cup\{e_{m+n-1}=v_{m+n-1}v_m\}$ so that $v_i\in V(P_m)$ for $i\in [1,m]$ and $v_j\in V(C_n)$ for $j\in [m,m+n-1]$. Note that $v_m$ is the vertex of degree 3. For $1\le i\le m+n-1$, define an edge labeling $f : E(T(m,n)) \to [1, m+n-1]$ by \[f(e_i) = \begin{cases}
m+n-(i+1)/2 & \mbox{ for odd } i,\\
i/2 & \mbox{ for even } i.
\end{cases}\]
We now have $$f^+(v_m) = \begin{cases}\frac{3(m+n)}{2} &\mbox{ for even } m,n,\\ \frac{3(m+n)-2}{2} &\mbox{ for odd } m,n, \\ \frac{3(m+n)-3}{2} &\mbox{ for odd } m \mbox{ and even } n, \\ \frac{3(m+n)-1}{2} &\mbox{ for even } m \mbox{ and odd } n. \end{cases}$$ Moreover, for $i\not= m$, $f^+(v_i) = m+n-1$ for odd $i$, $f^+(v_i)=m+n$ for even $i$. Thus, $\chi_{la}(T(m,n)\le 3$.

 Suppose there exists a local antimagic labeling $f$ that induces a 2-coloring of $T(m,n)$ with colors $x$ and $y$ such that $x<y$. Then $T(m,n)$ is bipartite so that $n$ is even. Let $X$ and $Y$ be the sets of vertices with colors $x$ and $y$, respectively. Clearly $||X|-|Y||\le 1$. Combining with Lemma~\ref{lem-2part}, we have $x|X|=(m+n)(m+n-1)/2=y|Y|$ and $|X|=|Y|+1$. By Lemma~\ref{lem-pendant}, $f^+(v_1)=x=m+n-1$. So $|X|=(m+n)/2$ and $|Y|=(m+n)/2-1$. Thus $y=(m+n)(m+n-1)/(m+n-2)$ which is not an integer, a contradiction. Thus, $\chi_{la}(T(m,n))= 3$.
\end{proof}

 For $a_1\ge a_2\ge \cdots\ge a_r\ge 3$, let $GB(a_1,a_2,\ldots,a_r)$ denote the generalized book graph which is the edge-gluing of cycles of order $a_i, 1\le i\le r$, at a common edge. We shall denote this common edge by $uv$ in the following three results.

\begin{lemma} For $r\ge 2$, $\chi_{la}(GB(a_1,a_2,\ldots,a_r))\ge 3$.
\end{lemma}
\begin{proof} Let $G=GB(a_1,a_2,\ldots,a_r)$. Suppose $G$ contains an odd cycle, then $\chi_{la}(G)\ge \chi(G)=3$. Suppose $G$ is bipartite, then $G$ has the same size of parts. By the contrapositive of Lemma~\ref{lem-2part}, we know $\chi_{la}(G)\not=2$. Therefore, $\chi_{la}(G)\ge 3$.
\end{proof}

\begin{theorem}\label{thm-GB3r} Suppose $r\ge 2$, we have $\chi_{la}(GB(3^{[r]}))=3$. \end{theorem}

\begin{proof} Let $G=GB(3^{[r]})$ such that $V(G)=\{u,v\}\cup\{x_i : 1\le i\le r\}$ and $E(G)=\{uv\}\cup\{ux_i ;1\le i\le r\}\cup\{vx_i : 1\le i\le r\}$. Define a bijection $f: E(G) \to [1,2r+1]$ by
\begin{enumerate}[(i)]
  \item $f(ux_i)=i$ for $1\le i\le r$,
  \item $f(vx_i)=2r+1-i$ for $1\le i\le r$,
  \item  $f(uv)=2r+1$.
\end{enumerate}
It is easy to verify that $f^+(x_i)=2r+1$ for $1\le i\le r$, $f^+(u)=r(r+1)/2+2r+1$ and $f^+(v)=(r+1)(3r+2)/2$. Hence, $f$ is a local antimagic labeling that induces a 3-coloring so that $\chi_{la}(G)\le 3$. Since $\chi_{la}(G)\ge \chi(G)= 3$, we have $\chi_{la}(G)=3$.
\end{proof}

 If $GB(a_1,a_2,\ldots,a_r)\not= GB(3^{[r]})$, it is easy to get a local antimagic labeling that induces a 4-coloring.

\begin{conjecture} If $a_1\ge 4$, then $\chi_{la}(GB(a_1,a_2,\ldots,a_r))= 4$.\end{conjecture}

 Let $G(a_1,a_2,\ldots,a_r;m)$ be obtained by identifying the vertex $u$ of\break $GB(a_1,a_2,\ldots,a_r)$ with a vertex of $m\ge 1$ copies of $P_2$.

\begin{theorem} Let $G = G(3^{[r]};m)$, then $$\chi_{la}(G)=\begin{cases}3 &\mbox{ if } G=G(3^{[r]};1) \mbox{ or } G(3^{[2]};2),\\ 4 &\mbox{ if } G=G(3^{[r]};2), r\ge 3, \\ m+1 &\mbox{ if }  m\ge  \binom{r}{2}\ge 3, \\ m+2 &\mbox{ if } 3\le m < \binom{r}{2}. \end{cases}$$ \end{theorem}

\begin{proof} For non-pendant vertices, we adopt the notations of $GB(3^{[r]})$. The pendant vertices are denoted by $y_j, 1\le j\le m$. By Theorem~\ref{thm-pendant}, we know $\chi_{la}(G)\ge m+1$. Since $G$ contains an odd cycle, we also have $\chi_{la}(G)\ge 3$. Suppose $m=1$. Define a bijection $f: E(G) \to [1,2r+2]$ by
\begin{enumerate}[(i)]
  \item $f(ux_i)=i$ for $1\le i\le r$,
  \item $f(vx_i)=2r+1-i$ for $1\le i\le r$,
  \item $f(uy_1)=2r+1$,
  \item $f(uv)=2r+2$.
\end{enumerate}

 Clearly, $f^+(u)=(r^2+9r+6)/2 \not= f^+(v)=(3r^2+5r+4)/2 \not= f^+(x_i)=2r+1 = f^+(y_1)$. Thus, $f$ is a local antimagic labeling that induces a 3-coloring. Thus, $\chi_{la}(G(3^{[r]};1)) = 3$.

 Consider $m=2$. Suppose $f$ is a local antimagic labeling that induces a 3-coloring. Without loss of generality, we must have $(r+1)(r+2)/2\le f^+(v)=f^+(y_1)\le 2r+3$. Hence,  $r=2$. The labeling $f(uv)=1, f(vx_1)=2, f(vx_2)=3, f(ux_1)=5, f(ux_2)=4, f(uy_1)=6, f(uy_2)=7$ gives $\chi_{la}(G(3^{[2]};2))=3$. For $r\ge 3$, we then have $\chi_{la}(G)\ge 4$. Define a bijection $f: E(G) \to [1,2r+3]$ by
\begin{enumerate}[(i)]
  \item $f(ux_i)=2r+2-i$ for $1\le i\le r$,
  \item  $f(vx_i)=i$ for $1\le i\le r$,
  \item $f(uy_j)=2r+1+j$ for $1\le j\le 2$,
  \item $f(uv)=r+1$.
\end{enumerate}

 Clearly, $f$ is a local antimagic labeling that induces a 4-coloring with $f^+(u)=(r+2)(3r+5)/2, f^+(v)=\binom{r+2}{2}$, $f^+(y_j)=2r+1+j, j=1,2$ and $f^+(x_i)=2r+2, 1\le i\le r$. Thus, $\chi_{la}(G(3^{[r]};2))=4$ for $r\ge 3$.

 Consider $m\ge 3$. We have $\chi_{la}(G)\ge m+1\ge 4$. Suppose $m\ge \binom{r}{2}$. Define a bijection $f: E(G) \to [1,2r+m+1]$ by
\begin{enumerate}[(i)]
  \item $f(ux_i)=2r+2-i$ for $1\le i\le r$,
  \item $f(vx_i)=i$ for $1\le i\le r$,
  \item $f(uy_j)=2r+1+j$ for $1\le j\le m$,
  \item $f(uv)=r+1$.
\end{enumerate}

 Clearly, $f$ is a local antimagic labeling that induces an $(m+1)$-coloring with $f^+(u)=(m+r+1)(3r+m+2)/2$, $f^+(v)=\binom{r+2}{2}=f^+(y_{\binom{r}{2}})$, $f^+(y_j)=2r+1+j, j\in [2,m]\backslash\{\binom{r}{2}\}$ and $f^+(x_i)=2r+2=f^+(y_1), 1\le i\le r$. Thus, $\chi_{la}(G(3^{[r]};m)) = m+1$ if $m\ge \binom{r}{2}\ge 3$.

 Suppose $m < \binom{r}{2}$, then $r\ge 4$. If $\chi_{la}(G)=m+1$, we may assume that $G$ admits a local antimagic labeling $f$ with $f^+(y_1)= f^+(x_i) \ge r(2r+1)/2, f^+(y_2)=f^+(v)\ge \binom{r+2}{2}$. Observe that $f^+(y_j)\le 2r+1+m < 2r+1+\binom{r}{2} = \binom{r+2}{2}\le f^+(v)$ for $1\le j\le m$, a contradiction. Thus, $\chi_{la}(G)\ge m+2$. Define a bijection $f: E(G)\to [1, 2r+m+1]$ by
\begin{enumerate}[(i)]
  \item $f(ux_i)=i+1$ for $1\le i\le r$,
  \item $f(vx_i)=2r+2-i$ for $1\le i\le r$,
  \item $f(uy_j)=2r+1+j$ for $1\le j\le m$,
  \item $f(uv)=1$.
\end{enumerate}

 It is easy to show that $f$ is a local antimagic labeling that induces a 4-coloring with $f^+(u)=(r+1)(r+2)/2 + m(4r+m+3)/2$, $f^+(v)=3r(r+1)/2 + 1$, $f^+(x_i)=f^+(y_2) = 2r+3, f^+(y_j)=2r+1+j$ for $j=1,3,4,\ldots, m$. Thus, $\chi_{la}(G(3^{[r]};m))=m+2$ if $3\le m< \binom{r}{2}$.
\end{proof}

\begin{problem} Study $\chi_{la}(G(a_1,a_2,\ldots,a_r;m))$ for $a_1\ge 4$. \end{problem}

 Suppose $G$ is of order $m$. Let $G\odot H$ be the graph obtained from $G$ and $m$ copies of $H$ by joining the $i$-th vertex of $G$ to each vertex of the $i$-th copy of $H$.

 Let $G=C_m \odot O_n$ with $V(G)=\bigcup\limits_{i=1}^m(\{v_{i,j} : 1\le j\le n\}\cup \{u_i\})$
and $E(G) =\bigcup\limits_{i=1}^m (\{u_iv_{i,j} : 1\le j\le n\}\cup\{e_i\})$, where $e_i = u_iu_{i+1}$ for $1\le i\le m$, and $u_{m+1}=u_1$ by convention. We shall keep these notation in the following discussion.

\begin{lemma}\label{lem-lowerCmOn} For $m\ge 3$ and $n\ge 1$, $\chi_{la}(C_m \odot O_n) \ge mn+2$.  \end{lemma}

\begin{proof} Let $f$ be a local antimagic labeling of $G=C_m \odot O_n$. Let $e$ be an edge of $G$ such that $f(e)=m(n+1)$ which is the size of $G$. If $e$ is not a pendant edge, then by Lemma~\ref{lem-pendant}, $c(f)\ge mn+2$. So we only need to deal with
$e=u_iv_{i,j}$ for some $i\in [1, m]$ and $j\in [1,n]$.
By renumbering we may assume that $e=u_1v_{1,1}$. Note that $f^+(v_{i,j})\le m(n+1)$.

 Suppose $f^+(u_i)$ and $f^+(u_{i+1})$ are greater than $m(n+1)$ for some $i$, $1\le i\le m$. Since they are distinct, $c(f)\ge mn+2$. So we may assume that the induced colors of any two consecutive vertices of $C_m$ do not both greater than $m(n+1)$.
Let $k$ be the number of vertices in $C_m$ whose induced color is less than or equal to $m(n+1)$. Thus, $m-1\ge k\ge \lceil m/2\rceil$. All edges in the cycle $C_{m}$ are incident to at least one of these $k$ vertices. So there are exactly $m+kn$ distinct edges incident to these $k$ vertices.

 Now $(m+kn)(m+kn+1)/2\le m(n+1)k$. Since $m\ge k+1$,
$k(n+1)(m+kn+1)<[k(n+1)+1](m+kn+1)=(k+1+kn)(m+kn+1)\le 2m(n+1)k$. Hence $m+kn+1<2m$ or $kn+1<m$. Since $k\ge \lceil m/2\rceil$, $n=1$. For this case, we have $(m+k)(m+k+1)\le 4mk$. This implies that $(m-k)^2+m+k\le 0$ which is impossible.
\end{proof}

\begin{theorem}\label{thm-CmOn} For $m\ge 3$ and $n\ge 1$, $\chi_{la}(C_m \odot O_n)= mn+2$ if both $m$ and $n$ are even; otherwise $mn+2\le \chi_{la}(C_m \odot O_n)\le mn+3$.   \end{theorem}

\begin{proof} Suppose $G=C_m\odot O_n$. By Lemma~\ref{lem-lowerCmOn}, we know $\chi_{la}(G)\ge mn+2$.

Consider $m=2h\ge 4$ and $n=2k\ge 2$. Define $f:E(G) \to [1,2h(2k+1)]$ by $f(e_i) = i$ for $1\le i\le 2h$ and

\begin{alignat*}{2}
f(u_1v_{1,1}) & = 4h+1, &  \quad f(u_1v_{1,2})& =6h, \\
f(u_{2i-1}v_{2i-1,1}) & =6h+3-2i,   & f(u_{2i-1}v_{2i-1,2})& =6h+2-2i, \quad 2\le i\le h;\\
f(u_{2i}v_{2i,1}) & =4h+1-2i,   & f(u_{2i}v_{2i,2})& =4h+2-2i, \quad 1\le i\le h;\\
f(u_{r}v_{r,2j-1}) & =2h(2j-1)+r, &    & 1\le r\le 2h, 2\le j\le k;\\
f(u_{r}v_{r,2j})& =2h(2j+1)+1-r, & & 1\le r\le 2h, 2\le j\le k.
\end{alignat*}

It is easy to check that $f$ is a bijection. It is also easy to verify that all pendants have different colors from $2h+1$ to $2h(2k+1)$.

 Now, for $2\le i\le h$, we have
\begin{align*}
f^+(u_{2i-1})& =f(e_{2i-1})+f(e_{2i-2})+f(u_{2i-1}v_{2i-1,1})+f(u_{2i-1}v_{2i-1,2})\\&\quad +\sum_{j=2}^k f(u_{2i-1}v_{2i-1,2j-1})+\sum_{j=2}^k f(u_{2i-1}v_{2i-1,2j})\\
&=
(2i-1)+(2i-2)+(6h+3-2i)+(6h+2-2i)\\&\quad +\sum_{j=2}^k [2h(2j-1)+(2i-1)]+\sum_{j=2}^k [2h(2j+1)+1-(2i-1)]\\ & = 12h+2+\sum_{j=2}^k [8hj+1]=12h+2+(k-1)(4hk+8h+1).
\end{align*}
We can also get that $f^+(u_{1}) =12h+2+(k-1)(4hk+8h+1)$ and $f^+(u_{2i})=8h+2+(k-1)(4hk+8h+1)$ for $1\le i\le h$. So we have $\chi_{la}(G)\le mn+2$. Hence, $\chi_{la}(G)= mn+2$.

Consider odd $m,n$. Let $A=(a_{i,j})$ be a magic $(m,n)$ rectangle involving the integers $[1, mn]$ (for the existence of magic rectangle, please see~\cite{Hagedorn}). Let $g$ be a local antimagic labeling of $C_m$ with $c(g)=3$. Now we define a labeling $f$ for $G$ by
\begin{center}
$f(e)=g(e)+mn$ for $e$ is an edge of $C_m$;\\
$f(u_iv_{i,j})=a_{i,j}$, for $1\le i\le m$ and $1\le j\le n$.\\
\end{center}
 Clearly $c(f)=mn+3$. So $\chi_{la}(G)\le mn+3$.

 Consider $m=2h+1\ge 3$ and $n=2k\ge 2$. Let $\alpha$ be the $k$-vector whose $i$-th coordinate is $i$, i.e., $\alpha=(1,2,\dots, k)$. Let $J$ be the $k$-vector whose coordinates are 1.
Let $B$ be a $(2h+1)\times 2k$ matrix whose $i$-th row is $((i-1)kJ +\alpha \ \ (4h+2-i)kJ+\alpha)$. Hence each row sum of $B$ is $(4h+1)k^2 + k(k+1) = 4hk^2 + 2k^2 + k$. Note that $B$ contains all integers in $[1,4hk+2k]$. Similar to the case of odd $m,n$, we will obtain a local antimagic labeling $f$ of $G$ with $c(f)=mn+3$. So $\chi_{la}(G)\le mn+3$.

 Consider even $m=2h\ge 4$ and odd $n=2k-1\ge 1$. Let $X=\{jk : 1\le j\le 2h\}$. We define a labeling $\phi : E(C_{2h})\to X$ by $\phi(e_{2i-1})= ik$ and  $\phi(e_{2i})= (h+i)k$, $1\le i\le h$. Then $\phi^+(u_j)=(h+j)k$ for $2\le j\le 2h$ and $\phi^+(u_{1})=(2h+1)k$.

 Let $C$ be a $(2h)\times (2k)$ matrix whose $i$-th row is $((4h-i)kJ+\alpha \ \ (i-1)kJ +\alpha)$. Hence each row sum of $C$ is $N=(4h-1)k^2+k(k+1)=4hk^2+k$. Let $C'=(c_{i,j})$ be a $(2h)\times (2k-1)$ matrix obtained from $C$ by deleting the last column of $C$. So the $i$-th row sum of $C'$ is $N-ik$, $1\le i \le 2h$. Now we shall label the pendant edges $u_i v_{i,j}$ by entries of a suitable row of $C'$, $1\le i\le 2h$ and $1\le j\le 2k-1$.

 Let $\psi: E(G)\setminus E(C_{2h})\to [1, 4hk]\setminus X$ defined by
$\psi(u_{2i} v_{2i,j})= c_{2i, j}$ for $1\le i\le h$; $\psi(u_{2i+1} v_{2i+1,j})= c_{2i-1,j}$, for $1\le i\le h-1$; and $\psi(u_1v_{1,j})=c_{2h-1, j}$, where $1\le j\le 2k-1$. Now $\psi^+(u_1)=N-(2h-1)k$; $\psi^+(u_{2i})=N-2ik$ for $1\le i\le h$ and $\psi^+(u_{2i+1})=N-(2i-1)k$ for $1\le i\le h-1$. Note that $\psi(u_3v_{3,k})=4hk$ which is the largest label.

 Let $f$ be the labeling of $G$ obtained by combining $\phi$ and $\psi$. Hence\\ $f^+(u_1)=\phi^+(u_1)+\psi^+(u_1)=(2h+1)k+[N-(2h-1)k]=N+2k$. \\ $f^+(u_{2i})=\phi^+(u_{2i})+\psi^+(u_{2i})=(h+2i)k+[N-2ik]=N+hk$ for $1\le i\le h$.\\
$f^+(u_{2i+1})=\phi^+(u_{2i+1})+\psi^+(u_{2i+1})=(h+2i+1)k+[N-(2i-1)k]=N+(h+2)k$ for $1\le i\le h-1$.

 Here $c(f)=2h(2k-1)+3$ if $h\ne 2$.

 For $h=2$, we redefine the labeling $\phi$ by $\phi(e_1)=k$, $\phi(e_2)=3k$, $\phi(e_3)=4k$ and $\phi(e_4)=2k$. Then $\phi^+(u_1)=3k$, $\phi^+(u_2)=4k$, $\phi^+(u_3)=7k$ and $\phi^+(u_4)=6k$. Hence  $f^+(u_1)=N$, $f^+(u_2)=N+2k$, $f^+(u_3)=N+6k$ and $f^+(u_4)=N+2k$. Here $c(f)=2h(2k-1)+3$ if $h= 2$.

 This completes the proof. \end{proof}

\begin{example}
Consider $G=C_3\odot O_1$. Denote $v_i=v_{i,1}$. Define $f(u_1u_2)=2$, $f(u_2u_3)=3$, $f(u_3u_1)=4$, $f(u_1v_1)=5$,  $f(u_2v_2)=1$ and $f(u_3v_3)=6$. Then the colors of vertices are 1,5,6,11,13. So $\chi_{la}(G)\le 5$. By Theorem~\ref{thm-CmOn} we have $\chi_{la}(C_3\odot O_1)= 5$. \rsq
\end{example}

\begin{example}According to the proof of Theorem~\ref{thm-CmOn} we have the labelings
for $C_4\odot O_1$ and $C_4\odot O_3$, respectively.\\
\centerline{\raisebox{5mm}[7mm][7mm]{\epsfig{file=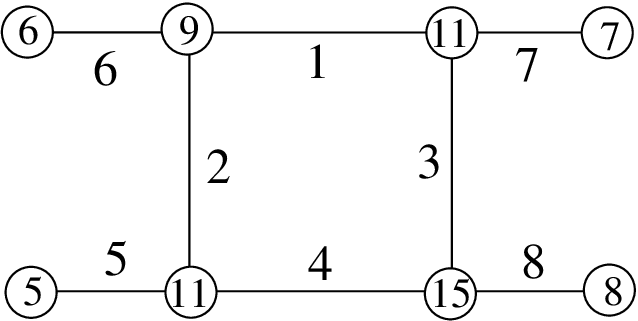, width=3.5cm}}\qquad \epsfig{file=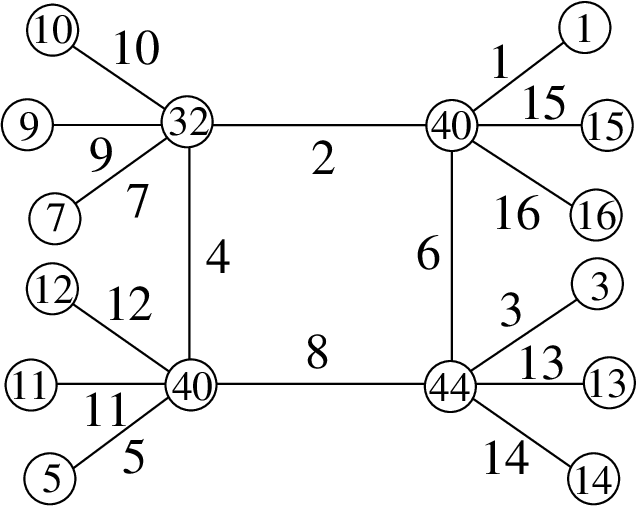, width=3.5cm}}

 After swapping some labels we have \\
\centerline{\raisebox{5mm}[7mm][7mm]{\epsfig{file=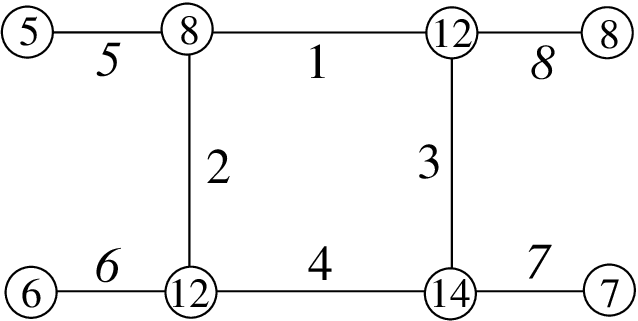, width=3.5cm}}\qquad \epsfig{file=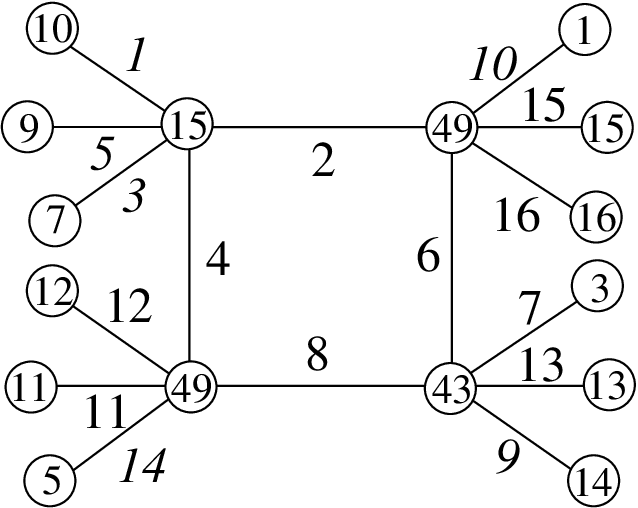, width=3.5cm}}
 So $\chi_{la}(C_4\odot O_1)\le 6$ and $\chi_{la}(C_4\odot O_3)\le 14$. By Theorem~\ref{thm-CmOn} we have $\chi_{la}(C_4\odot O_1)=6$ and $\chi_{la}(C_4\odot O_3)=14$.
\rsq \end{example}

 We are only aware, after obtaining the above theorem, that Arumugam et al.~\cite[Section 3]{Arumugam+L+P+W} have also obtained partial solutions on $\chi_{la}(C_m\odot O_n)$. Particularly, their Lemmas 3.6, 3.12 and 3.13 imply that $\chi_{la}(C_m\odot O_n)= mn+2$ for even $m\ge 4$ and $n\ge 3$. Moreover, Lemmas 3.14, 3.15 and 3.16 partially solved the case when $m$ is odd. This left $\chi_{la}(C_m\odot O_n)$ still unsolved for odd $m$ and finitely many $n$.


\section{$\chi_{la}(G)=|V(G)|$}

 For $m\ge 1$, $n_i\ge n_{i+1}$ $(1\le i\le m-1)$, and $n_1 + n_2 + \cdots + n_m \ge 1$, let $K(m;n_1,n_2,\ldots, n_m)$ be obtained from $K_m$ by joining $n_i$ pendant vertices to the $i$-th vertex of $K_m$. Note that $K(2;1,0)\cong P_3$ with $\chi_{la}(P_3)=3$ and $K(2;1,1)\cong P_4$ with $\chi_{la}(P_4)=3$ (see~\cite[Theorem 2.7]{Arumugam}). Moreover, $\chi_{la}(K(2;2,1))=4$ (see~\cite[Theorem 8]{LNS}). Observe that $K(1;n-1)\cong K(2;n-2,0)$ is the star graph $K_{1,n-1}$ of order $n$ with $\chi_{la}(K_{1,n-1})=n$. \\

\begin{theorem}~\label{thm-Km+pendant} For $m\ge 2$, $\chi_{la}(K(2;n_1,0))=n_1+2$.\\ Otherwise, $\chi_{la}(K(m;n_1,n_2,\ldots, n_m)) \le n_1 + n_2 + \cdots + n_m + m$ and the equality holds if and only if $(n_m+m-1)(n_m+m)/2 > n_1 + n_2 + \cdots + n_m + \binom{m}{2}$.
\end{theorem}

\begin{proof} Note that $G=K(m;n_1,n_2,\ldots, n_m)$ has order $n=n_1+n_2+\cdots + n_m + m$ and size $q=n_1+n_2+\cdots + n_m + \binom{m}{2}$.  By definition, it is easy to get $\chi_{la}(K(2;n_1,0))=n_1+2$. We now assume $G\not\cong K(2;n_1,0)$. Let  $V(G)=\{u_i : 1\le i\le m\}\cup\{u_{i,k}: 1\le i\le m, 1\le k\le n_i\}$ and  $E(G)=\{u_iu_j: 1\le i<j\le m\}\cup\{e_{i,k}=u_iu_{i,k}:1\le i \le m, 1\le k\le n_i\}$. Suppose $f$ is any local antimagic labeling of $G$. By definition, we must have all pendant vertex labels and all non-pendant vertex labels are mutually distinct respectively. Moreover, $f^+(u_i)\ne f^+(u_{i,k})$ and $f^+(u_{i,k}) \le q$ for all $i,k$. Since for $1\le i\le m$, $\deg(u_i)\ge n_m+m-1$, we also have $f^+(u_i) \ge 1+2+\cdots + (n_m+m-1) = (n_m+m-1)(n_m+m)/2$. Thus, if $(n_m+m-1)(n_m+m)/2 > q$, we have $f^+(u_i) > f^+(u_{j,k})$ for all $1\le i,j\le m, 1\le k\le n_j$. Therefore, $\chi_{la}=n$. This prove the sufficiency.

 To prove the necessity, suffice to show that if $(n_m+m-1)(n_m+m)/2 \le q$, then $c(f) < n$ for some local antimagic labeling $f$. Consider the sequence of edges $S=u_mu_{m-1}$, $u_mu_{m-2}$, $u_mu_{m-3}, \ldots$, $u_mu_1$, $e_{m,1}$, $e_{m,2}$, $e_{m,3}, \ldots$, $e_{m,n_m},$ $u_{m-1}u_{m-2}$, $u_{m-1}u_{m-3}, \ldots$, $u_{m-1}u_1$, $e_{m-1,1}$, $e_{m-1,2}$, $e_{m-1,3}, \ldots$, $e_{m-1, n_{m-1}}, \ldots$, $u_2u_1$, $e_{2,1}$, $e_{2,2}$, $e_{2,3}, \ldots$, $e_{2,n_2}$, $e_{1,1}$, $e_{1,2}, \ldots$, $e_{1,n_1}$. Define $f: S \to [1, q]$ according to the order in $S$. Observe that

\begin{enumerate}[(i)]
  \item $f$ is bijective,
  \item all the pendant vertex induced labels are distinct,
  \item for $1 \le i\le m, 1\le k\le n_i$, $f^+(u_i) > f^+(u_{i,k})$ if $n_i>0$,
  \item for $1\le i \le m-1$, $f(u_iu_j) > f(u_{i+1}u_j)$  $(j \ne i, i+1)$ and $\sum^{n_i}_{k=1}f(e_{i,k}) > \sum^{n_{i+1}}_{k=1}f(e_{i+1,k})$ so that $f^+(u_i) > f^+(u_{i+1})$.
\end{enumerate}

Thus, $f$ is a local antimagic labeling. We now have $f^+(u_m) = 1 + 2 + \cdots + (n_m + m - 1)  = (n_m+m-1)(n_m+m)/2$ and $f^+(u_{i,k})\le q$ for $1\le i\le m, 1\le k\le n_i$. Since $(n_m+m-1)(n_m+m)/2 \le q$, then there exists an edge $e$ with $f(e)=(n_m+m-1)(n_m+m)/2$, where $f$ is defined above. If $e=e_{i,k}$, then we have $f^+(u_{i,k})=f^+(u_m)$ so that $c(f) < n$. Otherwise, since for $1\le t\le m-1$, $f(u_tu_m)=m-t < m(m-1)/2 < (n_m+m-1)(n_m+m)/2=f(e)$, we must have $e=u_iu_j$ for $1\le i < j\le m-1$. We have the following two cases.

 {\bf Case (a).} $n_2=0$. In this case, $m\ge 3$. Note that $K_m$ has size $\binom{m}{2}=1+2+\cdots + (m-1)=f^+(u_m)$. Thus, $e=u_2u_1$ and $e_{1,1}$ must be the next unlabeled pendant edge. We now swap the labels of $u_2u_1$ and $e_{1,1}$. It is easy to verify that a new local antimagic labeling $g$ with $g^+(u_m)=g^+(u_{1,1})=\binom{m}{2}$ is obtained. Therefore, $c(g) < n_1 + m = n$.

 {\bf Case (b).} $n_2\ne 0$. In this case, according to our labeling sequence, $e_{j,1}$ must be the next unlabeled pendant edge.  Let $S'$ be obtained from $S$ by putting $e_{j,1}$ right before $e$. Now, define $g : S' \to [1,q]$ according to the order in $S'$. One can verify that all the observations under $f$ still hold under $g$. Moreover, $g^+(u_{j,1})=g^+(u_m)$. Thus, $g$ is a local antimagic labeling with $c(g) < n$. \end{proof}

 By Theorem~\ref{thm-pendant}, we know $\chi_{la}(K(2;a,b))\ge a+b+1$.

\begin{corollary}\label{cor-K2ab} For $a\ge b\ge 2$, \[\chi_{la}(K(2;a,b))=\begin{cases}
a+b+2 & \mbox{ if } a < b(b+1)/2 \\
a+b+1 & \mbox{ otherwise}.
\end{cases}\] \end{corollary}

 Suppose $k\ge2$, $n_1, \ldots, n_k\ge 0$ and $n_1 + \cdots + n_k\ge 2$. Let $Ct(k;n_1,\ldots,n_k)$ be the caterpillar graph obtained from the path $P_k=v_1v_2\cdots v_k$ by joining $n_i$ pendants to $v_i$. Consider the following two conditions:

\begin{enumerate}[{$C_1:$}]
  \item $\min\{(n_1+1)(n_1+2)/2, (n_2+2)(n_2+3)/2, (n_3+1)(n_3+2)/2\} > n_1+n_2+n_3+2$.
  \item No $n_1+n_3+2$ of distinct integers in $[1,n_1+n_2+n_3+2]$ can have sum of the $n_1+1$ integers equal sum of the remaining $n_3+1$ integers.
\end{enumerate}

\begin{theorem}\label{thm-Ct3} Suppose $n_1, n_2, n_3\ge 1$. If $Ct(3;n_1,n_2,n_3)$ satisfies conditions $C_1$ and $C_2$, then $\chi_{la}(Ct(3;n_1,n_2,n_3)) = n_1+n_2+n_3+3$. \end{theorem}

\begin{proof} Let $G=Ct(3;n_1,n_2,n_3)$ be the caterpillar graph obtained from the path $P_3=xyz$ by joining pendants $x_1, \dots, x_{n_1}$ to $x$, pendants $y_1, \dots, y_{n_2}$ to $y$ and pendants $z_1, \dots, z_{n_3}$ to $z$.
Let $f$ be a local antimagic labeling of $G$. Note that $f^+(x)\ge \frac{1}{2}(n_1+1)(n_1+2)$, $f^+(y)\ge \frac{1}{2}(n_2+2)(n_2+3)$ and $f^+(z)\ge \frac{1}{2}(n_3+1)(n_3+2)$. Moreover, $f^+(y)\ne f^+(x)$ and $f^+(y)\ne f^+(z)$. By $C_2$, we obtain that $f^+(x)\ne f^+(z)$.
Combining the results above, by $C_1$ we have $\chi_{la}(G)\ge n_1+n_2+n_3+3$.

 We now give a labeling $f: E(G) \to [1,n_1+n_2+n_3+2]$. By symmetry, we only need to consider three possibilities. Suppose $n_1\le n_2\le n_3$, we label in the sequence $xx_1, xx_2, \ldots,$ $xx_{n_1}, xy, yy_1, yy_2, \ldots, yy_{n_2}, yz, zz_1, zz_2, \ldots, zz_{n_3}$. Clearly, $f^+(z)=(n_3+1)(2n_1+2n_2+4)/2 > f^+(y)=$ $(n_2+2)(2n_1+n_2+3)/2 > f^+(x)=(n_1+1)(n_1+2)/2$ which in turn greater than all the pendant vertex labels. Thus, $c(f)=n_1+n_2+n_3+3$. Suppose $n_1\le n_3 < n_2$, we label in the sequence  $xx_1, xx_2, \ldots,$ $xx_{n_1}, xy, zz_1, zz_2, \ldots,$ $zz_{n_3}, yz, yy_1, yy_2, \ldots, yy_{n_2}$. Similarly, we have $f^+(y) > f^+(z) > f^+(x)$ which in turn greater than all the pendant vertex labels. Thus, $c(f)=n_1+n_2+n_3+3$. Finally, suppose $n_2 < n_1\le n_3$, we label in the sequence $yy_1, yy_2, \ldots, yy_{n_2}, xy, xx_1, xx_2, \ldots, xx_{n_1}, zz_1, zz_2, \ldots, zz_{n_3}, yz$. Similarly, we have we have $f^+(z) > f^+(x) > f^+(y)$ which in turn greater than all the pendant vertex labels. Thus, $c(f)=n_1+n_2+n_3+3$.
\end{proof}

 We are not able to find a $Ct(3;n_1,n_2,n_3)$ that satisfies Conditions $C_1$ and $C_2$. Thus, we have the following problem.

\begin{problem} Prove the existence of $Ct(3;n_1,n_2,n_3)$ in Theorem~\ref{thm-Ct3}. \end{problem}

 Let us consider a special case $Ct(3;n_1,0,n_3)$ with $n_3\ge n_1$. In this case, $C_1$ does not hold. When $n_1=1$, $Ct(3;1,0,n_3)$ is a coconut graph. It is known that $\chi_{la}(Ct(3;1,0,n_3))=n_3+2$. But we shall also include this result in the following corollary.

\begin{corollary} Suppose condition $C_2$ does not hold and $1\le n_1\le n_3$, then $\chi_{la}(Ct(3;n_1,0,n_3))=n_1+n_3+1$.
\end{corollary}
\begin{proof} Let $Ct(3;n_1,0,n_3)$ be as defined in Theorem~\ref{thm-Ct3}.
Since $C_2$ does not hold, $[1, n_1+n_3+2]$ has a bipartition $(S_1,S_3)$ such that $|S_1|=n_1+1$, $|S_3|=n_3+1$ and the sum of numbers in $S_1$ equals $\frac{1}{4}(n_1+n_3+2)(n_1+n_3+3)$ (of course, $n_1+n_3\equiv 1,2\pmod{4}$).

 Suppose $1\in S_1$. Choose $a\in S_3$ arbitrary. Define $f:E(Ct(3;n_1,0,n_3)) \to [1,n_1+n_3+2]$ such that $f(xy)=1, f(yz)=a$, $\{f(xx_i)\;|\; 1\le i\le n_1\} = S_1\setminus\{1\}$ and $\{f(zz_j)\;|\;1\le j\le n_3\}=S_3\setminus\{a\}$. Now $f^+(x)=f^+(z)=\frac{1}{4}(n_1+n_3+2)(n_1+n_3+3)$, $f^+(y)=1+a$ which is equal to a label of a pendant vertex. So we obtain that $\chi_{la}(Ct(3;n_1,0,n_3))\le n_1+n_3+1$. Similarly for $1\in S_3$.

 By Theorem~\ref{thm-pendant} we have $\chi_{la}(Ct(3;n_1,0,n_3)) = n_1 + n_3 + 1$.
\end{proof}

\begin{corollary} Suppose condition $C_2$ holds and $n_1\le n_3 < (n_1+2)(n_1-1)/2$, then $\chi_{la}(Ct(3;n_1,0,n_3))=n_1+n_3+2$.
\end{corollary}

\begin{proof} Under the assumption, $n_1> 2$. Let $Ct(3;n_1,0,n_3)$ be as defined in Theorem~\ref{thm-Ct3}. Condition $C_2$ holds implies that $f^+(x)\ne f^+(z)$ for all possible local antimagic labeling $f$ of $Ct(3;n_1,0,n_3)$. Moreover, $n_3 < (n_1+2)(n_1-1)/2$ implies that $(n_3+1)(n_3+2)/2 \ge (n_1+1)(n_1+2)/2 > n_1+n_3+2$ so that $f^+(x)$ and $f^+(z)$ are larger than all other pendant vertex colors for all possible local antimagic labeling of $Ct(3;n_1,0,n_3)$. Thus, $c(f)\ge n_1+n_3+2$.

 Define $f:E(Ct(3;n_1,0,n_3)) \to [1,n_1+n_3+2]$ such that $f(xy)=1, f(yz)=2$, $f(xx_i) = i+2$ for $1\le i\le n_1$, $f(zz_j) = n_1+2+j$ for $1\le j\le n_3$. Clearly, $f$ is a local antimagic labeling with $n_1+n_3+2$ distinct vertex colors.  Thus, $\chi_{la}(Ct(3;n_1,0,n_3)) = n_1 + n_3 + 2$.
\end{proof}

\begin{corollary} Suppose condition $C_2$ holds and $n_1\le  (n_1+2)(n_1-1)/2\le n_3 $, then $\chi_{la}(Ct(3;n_1,0,n_3))=n_1+n_3+1$.
\end{corollary}

\begin{proof} Under the assumption, $n_1\ge 2$. Define $f:E(Ct(3;n_1,0,n_3)) \to [1,n_1+n_3+2]$ such that $f(xy)=1$, $f(xx_i) = i+1$ for $1\le i\le n_1$, $f(yz)=n_1+2$ and $f(zz_j) = n_1+2+j$ for $1\le j\le n_3$. Now, $f^+(y)=n_1+3<f^+(x)=\frac{1}{2}(n_1+1)(n_1+2)$. Also $f^+(z)>f^+(y)$. Since $(n_1+2)(n_1-1)/2\le n_3$ is equivalent to $\frac{1}{2}(n_1+1)(n_1+2)\le n_1+n_3+2$, we have $f^+(x)=f^+(z_j)$ for some $j\ge 2$. Clearly, $f$ is a local antimagic labeling with $n_1+n_2+1$ distinct vertex colors.
By Theorem~\ref{thm-pendant} we have $\chi_{la}(Ct(3;n_1,0,n_3)) = n_1 + n_3 + 1$.
\end{proof}


\begin{problem} Study $\chi_{la}(Ct(k;n_1,\ldots,n_k))$. \end{problem}

 Note that the authors in~\cite{Arumugam+L+P+W} have also obtained results on $\chi_{la}(G\odot O_n)$ where $G$ is a path $P_m$ $(m\ge 2)$ or a complete graph $K_m$ $(m\ge 4)$. Note that $P_m\odot O_n \cong Ct(m;n^{[m]})$. In particular, they showed that $\chi_{la}(Ct(m;n^{[m]}))=|V(Ct(m;n^{[m]}))| = mn + m$ for $m\ge 3, n\ge 2$, and that $\chi_{la}(C_3\odot O_n) = |V(C_3\odot O_n)| = 3n+3$ for $n\ge 2$.

  Up to now, it is known that for a graph $G$ of order $n$, $\chi_{la}(G)=n$ if $G=K_n$ $(n\ge 2)$, $K_{1,n-1}$ $(n\ge 3)$, or $K(m;n_1,n_2,\ldots, n_m)$ for $(n_m+m-1)(n_m+m)/2 > n_1 + n_2 + \cdots + n_m + \binom{m}{2}$, or $Ct(3;n_1,n_2,n_3)$ of Theorem~\ref{thm-Ct3}, or $Ct(m;n^{[m]})$ $(m\ge 3, n\ge 2)$ or $C_3\odot O_n$ ($n\ge 2)$. It is also easy to verify that if $G$ is a graph of order $3 \le n\le 6$, then $\chi_{la}(G)=n$ if and only if $G=K_n, K_{1,n-1}$ or $K(2;2,2)$. In~\cite[Problem 3.1]{Arumugam}, the authors posed the problem: Characterize the class of graph $G$ of order $n$ for which $\chi_{la}(G)=n$. We end this paper with the following conjecture.

\begin{conjecture} A graph $G$ of order $n$ has $\chi_{la}(G)=n$ if and only if $G=K_n$ $(n\ge 3)$; or $K_{1,n-1}$ $(n\ge 3)$; or $K(m;n_1,n_2,\ldots, n_m)$ for $(n_m+m-1)(n_m+m)/2 > n_1 + n_2 + \cdots + n_m + \binom{m}{2}$; or $Ct(3;n_1,n_2,n_3)$ of Theorem~\ref{thm-Ct3}; or $Ct(m;n^{[m]})$ $(m\ge 3, n\ge 2)$; or $C_3\odot O_n$ $(n\ge 2)$. \end{conjecture}

\section*{Acknowledgments}
The authors wish to thank the referee for the valuable comments.

\end{document}